\documentclass[11pt]{scrartcl}
\usepackage[utf8]{inputenc}
\usepackage{microtype}
\usepackage{amsmath, amssymb, amsfonts, amsthm, mathtools}
\usepackage{MnSymbol}
\usepackage{graphicx}
\theoremstyle{definition}
\usepackage[backend=bibtex, %use bibtex as the backend not biber 
    style=numeric-comp,   % numeric-comp, reading for abstracts reference style
    sorting=nyt,          % sort Name, Year, Title
    citestyle=numeric-comp,    % A compact variant of numeric style
    maxnames=5,           % max names before et. al. in citation
    minnames=4,           % min names before et. al. in citation
    maxbibnames=99,       % forces print all author names in biblio
    hyperref=true,        % hyperlink citations
    doi=false,            % do not print doi
    isbn=false,           %do not print isbn
    url=false]             % do not print document url
    {biblatex}            % use biblatex
\usepackage{float}
\usepackage{subfloat}
\usepackage[config, labelfont={sf,bf},textfont=sf]{caption,subfig}  %Have multiple images side by side, each individually captioned.
    \DeclareCaptionTextFormat{bookmarktable}{\belowpdfbookmark{Table \thetable: #1}{\thetable}#1}
\usepackage[hidelinks]{hyperref} 
\hypersetup{colorlinks = true}
\newcommand{\footremember}[2]{%
\footnote{#2}
\newcounter{#1}
\setcounter{#1}{\value{footnote}}%
}

\author{
Bradley W. Brock\footremember{Brock} {CCR, La Jolla, CA, 92121},
 Robert Compton\footremember{Compton} 
{HC 1, Box 3024, Joshua Tree, CA 92252 Email: \url{rccompt@gmail.com}},
Warwick de Launey\footremember{Warrick} {CCR, La Jolla, CA, 92121}\\
and Jennifer Seberry\footremember{UoW} {School of Computing and Information Technology,  University of Wollongong, NSW 2522, Australia. Email: \url{jennifer_seberry@uow.edu.au}
}}
\title{On Generalized Hadamard Matrices and Difference Matrices: $Z_6$}

\begin{document}

\maketitle

\begin{abstract}
We give some very interesting matrices which are orthogonal over groups and, as far as we know, referenced, but in fact undocumented. This note is not intended to be published but available for archival reasons.

\end{abstract}


\section{Introduction}\label{sec:introduction}
(Caution: please note Jennifer Seberry has tried but been unable to contact the persons named as co-authors. Apolgogies for half finished results.) We note that the literature on this subject is very disorganized. Authors have not read the literature on their own key-words and papers thus ignored, on the other hand papers have been claimed to be published which have not. We have put together this note to make them accessible to all.

We do not give examples of all the literature but point to references which may help.

This work is compiled from journal work of the authors from old pieces of paper.

We first note that the matrices we study here have elements from groups, abelian and non-abelian, and may be written in additive or multiplicative notation. The matrices may have real elements, elements $\in \{1,-1\}$, elements $|n| \leq 1$, elements $\in \{1,i,~i^2=-1\}$, elements $\in \{1,i,-1,-i,~i^2=-1\}$, integer elements $\in \{a+ib,~i^2=-1\}$, $n$th roots of unity, the quaternions \{1 and $i,j,k,~i^2=j^2=k^2 =-1,ijk=-1\}$, $(a+ib) + j(c+id)$, \textit{a, b, c, d, integer} and \textit{i, j, k} quaternions or otherwise as specified.

We use the notations $B^{\top}$ for the transpose of $G$, $B^H$ for the group transpose, $B^C$ for the complwx conjugate of $B^{\top}$ and $B^V$ for the quaternion conjugate transpose.

In all of these matrices the inner product of distinct rows a and b is a - b or $a.b^{-1}$ depending on whether the group is written in additive or multiplicative form.

{\flushleft \bf LEXICON 1}\label{def:orthog}
\begin{enumerate}
\item \textbf{classical orthogonality:} 
$ H =\begin{bmatrix} 1& 1& 1& 1 \\ 1& 1& -1& -1\\ 1&-1&1&-1\\1&-1&-1&1\end{bmatrix}; \quad HH^{\top} =4I_4.$ 

\item \textbf{group orthogonality:}\[G = \begin{bmatrix}
1 & 1 & 1 & 1\\ 1 & a & b & ab\\ 1 & b  & ab & a \\ 1 & ab & a & b 
\end{bmatrix}; \qquad GG^H = (group)I_4 = (Z_2 \times Z_2)I\]
orthogonality is because the inner product of distinct rows is the whole group the same number of times. The inner product is $\{g_{i1}g_{j1}^{-1},\dots,g_{in}g_{jn}^{-1}\}$

$GG^H = 5Z_4I_{20}$
\[\begin{array}{cccc|cccc|cccc|cccc|cccc}
0& 0& 0& 0&  0& 0& 0& 0&  0& 0& 0& 0&  0& 0& 0& 0&  0& 0& 0& 0\\
0& 0& 0& 0&  2& 2& 2& 2&  4& 4& 4& 4&  1& 1& 1& 1&  3& 3& 3& 3\\
0& 0& 0& 0&  3& 3& 3& 3&  1& 1& 1& 1&  4& 4& 4& 4&  2& 2& 2& 2\\
0& 0& 0& 0&  4& 4& 4& 4&  3& 3& 3& 3&  2& 2& 2& 2&  1& 1& 1& 1\\ \hline
0& 2& 3& 4&  3& 4& 0& 1&  4& 0& 1& 2&  0& 1& 2& 3&  1& 2& 3& 4\\
0& 2& 3& 4&  4& 1& 3& 0&  0& 2& 4& 1&  1& 3& 0& 2&  2& 4& 1& 3\\
0& 2& 3& 4&  0& 3& 1& 4&  1& 4& 2& 0&  2& 0& 3& 1&  3& 1& 4& 2\\
0& 2& 3& 4&  1& 0& 4& 3&  2& 1& 0& 4&  3& 2& 1& 0&  4& 3& 2& 1\\ \hline
0& 4& 1& 3&  4& 0& 1& 2&  1& 3& 0& 2&  3& 1& 4& 2&  0& 4& 3& 2\\
0& 4& 1& 3&  0& 2& 4& 1&  0& 1& 2& 3&  4& 3& 2& 1&  4& 2& 0& 3\\
0& 4& 1& 3&  1& 4& 2& 0&  3& 2& 1& 0&  1& 2& 3& 4&  3& 0& 2& 4\\
0& 4& 1& 3&  2& 1& 0& 4&  2& 0& 3& 1&  2& 4& 1& 3&  2& 3& 4& 0\\ \hline
0& 1& 4& 2&  0& 1& 2& 3&  3& 1& 4& 2&  0& 4& 3& 2&  4& 1& 3& 0\\
0& 1& 4& 2&  1& 3& 0& 2&  4& 3& 2& 1&  4& 2& 0& 3&  3& 4& 0& 1\\
0& 1& 4& 2&  2& 0& 3& 1&  1& 2& 3& 4&  3& 0& 2& 4&  1& 0& 4& 3\\
0& 1& 4& 2&  3& 2& 1& 0&  2& 1& 4& 3&  2& 3& 4& 0&  0& 3& 1& 4\\ \hline
0& 3& 2& 1&  1& 2& 3& 4&  0& 4& 3& 2&  4& 1& 3& 0&  2& 0& 3& 1\\
0& 3& 2& 1&  2& 4& 1& 3&  4& 2& 0& 3&  3& 4& 0& 1&  3& 2& 1& 0\\
0& 3& 2& 1&  3& 1& 4& 2&  3& 0& 2& 4&  1& 0& 4& 3&  0& 1& 2& 3\\
0& 3& 2& 1&  4& 3& 2& 1&  2& 3& 4& 0&  0& 3& 1& 4&  1& 3& 0& 2
\end{array}
\]
Example of $GH(20,Z_4)$.

\item \textbf{Butson orthogonality:} \cite{AB1962, WCRN2007}
\[U= \begin{bmatrix} 1 & 1  & 1\\ 1& \omega & \omega^2 \\ 1 & \omega^2 & \omega
\end{bmatrix}; \quad BB^C =3I_3, \quad w^3 = 1, \quad 1+w+w^2 = 0\]
orthogonality depends on the fact that the $n$ $n$th roots of unity add to zero.

\item \textbf{complex orthogonality (A):} 
\[ A = \begin{bmatrix} 
1 & 1\\ i & -i
\end{bmatrix}; \quad BB^C =2I_2,\]
orthogonality is independent of the internal elements.

\item \textbf{complex orthogonality (B):}
\[B = \begin{bmatrix}
i & 1\\ 1  & i
\end{bmatrix}; \quad BB^H =2I_2,\]
$ B^H $is $B$ transposed complex conjugate.

\item \textbf{quaternion orthogonality:} \cite{SFAWW}
\[V = \begin{bmatrix}
 1& k \\ i& j \end{bmatrix}; \quad VV^{Q} =ZI_2, \quad i^Q = -i. \]
 $i,j,k$ are quaternions.
\end{enumerate}

{\flushleft \bf LEXICON 2}\label{trix} is an orthogonal matrix of order $n$ with entries 1, -1. \cite{WSW} 
\begin{enumerate}

\item 
A \textbf{Cretan matrix} \cite{JSNB2015b} is an orthogonal matrix with entries $\leq 1$. They are also referred to in the literature as Fermat, Hadamard, Mersenne, Euler, Belevitch and conference matrices.

\item 
A \textbf{difference matrix} \cite{DJ1980} of size $v \times b$ has entries from a group. The inner product ${\bf}a \cdot {\bf}b^{-1}$ of any pair of distinct rows is a constant number of each of the elements of the group.
\[\begin{bmatrix}
0& 0& 0& 0& 0& 0\\
0& 0& 1& 2& 2& 1\\
0& 1& 0& 1& 2& 2\\
0& 2& 1& 0& 1& 2\\
0& 2& 2& 1& 0& 1\\
0& 1& 2& 2& 1& 0
\end{bmatrix} \quad \text{is a }GH(6,Z_3) \text{ difference matrix}
\]

\item  A \textbf{generalized/generalized Hadamard matrix} \cite{AB1962, AB1963}
is $v \times v$ difference matrix.

\item  A \textbf{symmetric balanced incomplete block design: SBIBD}\cite{}
is the incidence matrix of a $(v,k,\lambda)$-design or $(v,k,\lambda)$-configuration (Ryser \cite{Ryser}): it is a
$v \times v$ matrix with entries 0, 1, where 1 occurs $k$ times in each row and column and the inner product of any distinct pairs of rows and columns is $\lambda$.

\item A \textbf{balanced incomplete block design: BIBD}\cite{JRSMY1992, WSW}
A BIBD, $(v,k,r, k, \lambda)$-design is $v \times b$ matrix with entries 0, 1, where 1 occurs $r$ times in each row and $k$ times in each column and the inner product of any distinct pairs of rows  is $\lambda$.

\item A \textbf{Bhaskar Rao Design}\cite{WD1989}
 design has entries from a group and zero so that the inner product ${\bf}a \cdot {\bf}b$ of any pair of distinct rows is a constant number of each of the elements of the group. The underlying design where the group elements are replaced by 1 is the incidence matrix of an BIBD or balanced incomplete block design.

\item \textbf{Butson-Hadamard Matrices} \cite{AB1962, AB1963} are generalized Hadamard  matrices where the group is the $n$th roots of unity.

\item \textbf{Quaternion orthogonal matrices} \cite{SFAWW} are generalized Hadamard  matrices where the group is the quaternions.
\end{enumerate}

\subsection{Other Constructions}

Warwick de Launey's research featured multiple articles on generalized Hadamard matrices \cite{WD1983,WD1984b,WD1986a,WD1987b,WD1989,WDL1992}, articles with Ed Dawson \cite{WDLED1992,WDLED1994}, Charles J. Colbourn \cite{CCWDL1996}, Rob Craigen \cite{RCWDL2009}. 

Other authors who have multiple articles include Jennifer Seberry \cite{JRS1979b,JRS1980a,JRS1986} and with Geramita \cite{AGJS1979}, Yamada \cite{JRSMY1992}, W.D. Wallis and Anne Street \cite{WSW}, and Balonin \cite{JSNB2015b}. Ferenc Sz\"{o}ll\H{o}si \cite{FS2010b,FS2011a,FS2012a,FS2012b} and with Lampio and \"{O}stergard \cite{PLFSPO2012}, Matolcsi \cite{MMFS2007}. Gennian Ge \cite{GG2005} and with Miao and Sun \cite{GGYM2010}. Dieter Jungnickel \cite{DJ1980} and with Grams \cite{DJGG1986}. 
 
For other constructions  and additional reading see: Ayd{\i}n and Altay \cite{CABA2013}, Camp and Nicoara \cite{WCRN2007}, Colbourn and Kreher \cite{CCDK1996}, Dawson \cite{JED1985}, Delsarte and Goethals \cite{PDJG1969}, Drake \cite{DD1979}, Egiazaryan \cite{KE1984}, Evans \cite{AE1987}, Hadamard \cite{JH1893}, Harada, Lam, Munemasa and Tonchev \cite{MHCL2010}, Hayden \cite{JH1997}, Hiramine \cite{YH2014} and with Suetake \cite{YHCS2013}, Horadam \cite{KH2000}, Karlsson \cite{BK2009,BK2011}, Lampio and \''{O}sterg\.{a}rd \cite{PLPO2011}, Liu \cite{ZL1977}, Sharma and Sookoo \cite{BSNS2010}, Shrikhande \cite{SS1964}, Deborah Street \cite{DS1979}, Todorov \cite{DT2012}, Winterhof \cite{AW2000}, and Zhang \cite{YZ2010} and with Duan, Lu, and Zheng \cite{YZLD2002}.

Another useful reference is \fullcite{La-Jolla}.

\section{A Summary of Known Generalized Hadamard  Matrices over \texorpdfstring{$GH(n;Z_6)$}{GH(n;Z6}}
\begin{table}[H]
\caption{Table of Existence of \texorpdfstring{$GH(n;Z_6)$ $2 \leq n \leq 52$}{GH(n;Z6 2 ≤ n ≤ 52)}}
\label{table:GH-Z_6}
\begin{tabular}{cl|cl|cl}
\hline $n$ & Comment             & $n$ & Comment              & $n$ & Comment \\ \hline
       2 & $H(2;Z_2)$ exists     & 19 & ?                     & 36 & $H(36;Z_2)$ exists \\ 
       3 & $GH(3; Z_3)$ exists   & 20 & $10 \times 2$         & 37 & ? \\ 
       4 & $H(4; Z_2)$ exists    & 21 & $GH(7; Z_3)$ (\ref{subS:GH(7;Z_6)})   & 38 & ? \\ 
       5 & NE                    & 22 & ?                     & 39 & ? \\ 
       6 & $GH(6; Z_3)$ exists   & 23 & NE                    & 40 & $10 \times 4$ \\ 
       7 & $GH(7; Z_6)$ (\ref{subS:GH(7;Z_6)})    & 24 & $H(24; Z_2)$ exists   & 41 & NE \\ 
       8 & $H(8; Z_2)$ exists    & 25 & ?                     & 42 & $7 \times 6$ \\ 
       9 & $3 \times 3$          & 26 & ?                     & 43 & ? \\ 
      10 & $GH(10; Z_3)$ (\ref{subS:GH(10;Z_6)})   & 27 & $3 \times 3 \times 3$ & 44 & $H(44;Z_2)$ exists \\ 
      11 & NE                    & 28 & $H(28; Z_2)$ exists   & 45 & NE \\ 
      12 & $H(12; Z_2)$ exists   & 29 & NE                    & 46 & ? \\ 
      13 & NE                    & 30 & $10 \times 3$         & 47 & NE \\ 
      14 & $7 \times 2$          & 31 & ?                     & 48 & $H(48; Z_2)$ exists \\ 
      15 & NE                    & 32 & $H(32; Z_2)$ exists   & 49 & $7 \times 7$ \\ 
      16 & $H(12; Z_2)$ exists   & 33 & NE                    & 50 & ? \\ 
      17 & NE                    & 34 & $GW(17,16,15;Z_3)$ exists & 51 & ? \\ 
      18 & $3 \times 3 \times 2$ & 35 & NE                    & 52 & $H(52; Z_2)$ exists \\ 
\hline 
\end{tabular}
\end{table}

\section{Some unpublished interesting matrices}\label{S:Interesting-Matrices}
\subsection{Matrix \texorpdfstring{$GH(7; Z_6)$}{GH(7;Z6)}}\label{subS:GH(7;Z_6)}

\begin{tabular}[c]{ccc}
$\left[\begin{array}{r|rrr|rrr}
-1 & 1        & 1         & 1       & 1        & 1        & 1\\ \hline
1  & -\omega  & \omega    & \omega  & \omega^2 & 1        & 1\\
1  & \omega   & -\omega   & \omega  & 1        & \omega^2 & 1\\
1  & \omega   & \omega    & -\omega & 1        & 1        & \omega^2 \\ \hline
1  & \omega^2 & 1         & 1       & -\omega  & \omega   & \omega \\
1  & 1        & \omega^2 & 1       & \omega   & -\omega  & \omega \\
1  & 1        & 1         & \omega^2& \omega   & \omega   & -\omega 
\end{array}\right]$ & and &
$\left[\begin{array}{r|rrr|rrr}
-1 & 1        & 1       & 1       & 1        & 1        & 1\\ \hline
1  & -\omega  & 1       & \omega  & \omega^2 & \omega   & 1\\
1  & 1        & -\omega & 1       & \omega   & \omega^2 & \omega\\
1  & \omega   & 1       & -\omega & 1        & \omega   & \omega^2\\ \hline
1  & \omega^2 & \omega  & 1       & -\omega  & 1        & \omega\\
1  & \omega   & \omega^2& \omega  & 1        & -\omega  & 1\\
1  & 1        & \omega  & \omega^2& \omega   & 1        & -\omega 
\end{array}\right]$
\end{tabular}

\subsection{Matrix \texorpdfstring{$GH(10;Z_6)$}{GH(10;Z6)}}\label{subS:GH(10;Z_6)}
Let $X$, $Y$, $Z$, $W$ be the $5 \times 5$ circulant matrices with first rows: 
\[\begin{matrix*}[r]
-1, & \omega,   & \omega^2, & \omega^2, & \omega\\
1,  & \omega^2, & \omega,   & \omega,   & \omega^2
\end{matrix*} \qquad \qquad
\begin{matrix*}[r]
1, & \omega,    & \omega^2,  & \omega^2, & \omega\\
1, & -\omega^2, & -\omega^2, & -\omega,  & -\omega^2 
\end{matrix*}
\]
Then
\[\begin{bmatrix}
X & Y\\ W & Z
\end{bmatrix}
\]
is a $GH(10;Z_6)$.

\subsection{Circulant matrices}\label{subS:circulant-matrices}
\begin{enumerate}
\item First rows of $GH(20; Z_4)$ from 4 circulant matrices
      \[\begin{matrix*}[r]ab & b  & e & e & e\end{matrix*}, \qquad
      \begin{matrix*}[r] b & a  & e & e & a\end{matrix*}, \qquad
      \begin{matrix*}[r]e  & ab & e & e & b\end{matrix*}, \qquad
      \begin{matrix*}[r]ab & a  & e & e & a\end{matrix*}.\]

\item Let $C = GW(17, 16; Z_3)$ then the following is a $GH(34; Z_6)$
      \[\begin{bmatrix*}[r]
       I+C & I-C\\ I-C^{\ast} & -I-C^{\ast}\end{bmatrix*}\]

\item First rows of $GH(28; Z_4)$ from 4 circulant matrices
%\marginpar{\textsf{Can we find for 3, 5, 9, etc.}}
      \[\begin{array}{lcl}
      \begin{matrix*}[r]e & a & a & ab & a & ab & ab\end{matrix*}\,, &&
      \begin{matrix*}[r]e & b & b & ab & b & ab & ab\end{matrix*}\,,\\
      \begin{matrix*}[r]e & a & a & b & a & b & b\end{matrix*}\,, &&
      \begin{matrix*}[r]e & ab & b & ab & a & a & b\end{matrix*}\,.
      \end{array}\]

%\item Notes on $D(p,q,s,t) = \text{circ}(e,q,s,q,t,t,s)$. All $D(e,q,s,t)D(e,q,s_{\ell},t)^{\top} = 2e+2a+2b+??$\marginpar{couldn't read end of line}
%      \[\begin{matrix*}[l]
%      D(e,a,b,ab) & D(ab,ab,b,a) & D(ab,ab,a,b) & D(b,b,a,ab)\\ 
%      D(e,a,b,ab) & D(b,b,ab,e)  & D(a,a,ab,e)  & D(a,a,b,e)\\ 
%      D(e,a,b,ab) & D(a,a,e,ab)  & D(e,e,b,a)   & D(ab,ab,e,b)\\
%      D(e,a,b,ab) & D(e,e,a,b)   & D(b,b,e,ab)  & D(e,e,ab,a)
%      \end{matrix*}\]

\item Brock's vectors of length 7 - these work but a  second example was corrupted. Gives $GH(21; Z_3)$.
      \[\begin{matrix*}[l]
       1121121  & 0012210  & 1012210\\
       0012210  & 0100001  & 2012210\\
       1012210  & 2012210  & 2220022
       \end{matrix*}\]
       
\item Brock's vectors of length 13 - How much works? Part of $GH(39; Z_3)$. \cite{AB1962,AB1963,WCRN2007,RCWDL2009}
      \[\begin{matrix*}[l]
       1200020020002 & 0011202202110 & 1011202202110\\
       0011202202110 & 0222121121222 & 2011202202110\\
       1011202202110 & 2011202202110 & 2100111111001
       \end{matrix*}\]

\item  Let $C_0$, $C_1$, $C_2$ be the cubic residues classes of 13 then with $e$, $\omega$, $\omega^2$ the elements of $Z_3$ the first row for 13 is
      \[eI + \omega^2C_0+ eC_1 + eC_2\,, \qquad eI + eC_0 + \omega C_1 + \omega^2C_2\,, \qquad \omega I + eC_0 + \omega C_1+ \omega^2C_2\,.\]
\end{enumerate}
%\marginpar{This is where the edited pages ended, so I must have left a page behind}

hence if $\{q,s,t\} = \{a,b,ab\}$ from $Z_2 \times Z_2$ all
\[D(e, q, s, t)D(e, q, s, t)^{\top} = 2e + 2a + 2b + 2ab\]

\textbf{Cubic residues of 13}
$C_0 = \{1, 5, 8,12\}$ $C_1 = \{2, 3, 10, 11\}$ $C_2 = \{4, 6, 7, 9\}$.

\begin{align*}
C_0 + C_0 &= C_1+2C_2+4I \\
C_0 +C_1 &= C_0 +2C_1+ C_2 = C_1 + C_O\\
C_0 + C_2 &= 2C_O + C_1 + C_2 = C_2 + C_O 
\end{align*}
So Brock's $3 \times 13$ matrix is
\[\begin{matrix*}[l]
\omega I + \omega^2C_0 + eC_1 + eC_2 & eI + eC_0 + \omega C_1 + \omega^2C_2 & \omega I + eC_0 +\omega^2C_2\\
eI + \omega^2C_0 + \omega^2C_1 + \omega^2C_2 & eI + eC_0 + \omega C_1 + \omega^2C_2 & \omega^2I + eC_0 + \omega C_1+ \omega^2C_2\\
\omega I + eC_0 + \omega C_1 + \omega^2C_2 & \omega^2I+ eC_0 + \omega C_1 + \omega^2C_2 & \omega^2I + \omega C_0 + \omega C_1 + \omega C_2
\end{matrix*}
\]

\textbf{Other Strange Squares}
\begin{enumerate}
\item
\[\text{circ}(*, 1, -, 1, 0, -, -, 1, 1, 1, 0, 1, 0) \]
is orthogonal if '*' is ignored.
\item \[\text{circ}(0, a, e, ab, 0, e, e, ab^2, a, ab^2, 0, ab, 0)\]
is the $GW(13,9,6; S_3)$.

\item Let $C = GW(17, 16; Z_3)$ then the following is a $GH(34; Z_6)$
\[\begin{bmatrix*}[c]
I + C & I - C\\
I - C^{\ast} & -I - C^{\ast}
\end{bmatrix*}\]

\end{enumerate}

\section{Acknowledgements}

 The authors wish to thank Professor Rob Craigen, University of Manitoba, for his valuable comments. The authors also wish to sincerely thank Mr Max Norden (BBMgt) (C.S.U.) for his work preparing the layout and LaTeX version.

%\bibliographystyle{unsrt}
%\end{document}
\printbibliography

\begin{thebibliography}{99}

\bibitem{CABA2013}
Cafer Aydın and Bilal Altay. Domain of generalized difference matrix B(r,s) on some Maddox’s spaces. In: \textit{Thai J. Math}. 11.1 (2013), pp. 87–102.

\bibitem{BB1988}
B.W. Brock. Hermitian congruence and the existence and completion of generalized Hadamard matrices. In: \textit{J. Combin. Theory} Series A 49 (1988), pp. 233–261.

\bibitem{AB1962}
A.T. Butson. Generalized Hadamard matrices In: \textit{Proc. Amer. Math. Soc.} 3 (1962), pp. 894–898.

\bibitem{AB1963}
A.T. Butson. Relations among generalized Hadamard matrices, relative difference sets, and maximal length recurring sequences. In: \textit{Canad. J. Math}. 15 (1963), pp. 42–48.

\bibitem{WCRN2007}
W. Camp and R. Nicoara. “Subfactors and Hadamard Matrices”. In: ArXiv e-prints (2007). arXiv: 0704.1128 [math.OA].

\bibitem{CCDK1996}
Charles J. Colbourn and Donald L. Kreher. Concerning difference matrices. In: \textit{Des. Codes Cryptogr.} 9.1 (1996). Second Upper Michigan Combinatorics Workshop on Designs, Codes and Geometries (Houghton, MI, 1994), pp. 61–70.

\bibitem{CCWDL1996}
Charles J. Colbourn and Warwick de Launey. Difference Matrices. In: \textit{The CRC Handbook of Combinatorial Design.} Ed. by C.J. Colbourn and J.H. Dinitz. Boca Raton, FL: CRC Press, 1996. Chap. IV.11, pp. 287–297.

\bibitem{RCWDL2009}
Rob Craigen and Warwick De Launey. Generalized Hadamard matrices whose transposes are not generalized Hadamard matrices. In: \textit{J. Combinatorial Designs} 17.6 (2009), pp. 456–458.

\bibitem{JED1985}
J.E. Dawson. A construction for generalized Hadamard matrices $GH(4q,EA(q))$. In: \textit{J. Statist. Plan. Inf.} 11 (1985), pp. 103–110.

\bibitem{PDJG1969}
P. Delsarte and J.M. Goethals. Tri-weight codes and generalized Hadamard matrices. In: \textit{Inf. and Contr.} (1969), pp. 196–206.

\bibitem{DD1979} 
D.A. Drake. Partial $\lambda$-geometries and generalised Hadamard matrices. In: \textit{Canadian J. Math}. 31 (1979), pp. 617–627.

\bibitem{KE1984}
K.O. Egiazaryan. $n$-dimensional generalized Hadamard matrices. In: \textit{Akad. Nauk Armyan. SSR Dokl}. 78.5 (1984). In Russian. Armenian summary., pp. 203–207.

\bibitem{AE1987}
Anthony B. Evans. Difference matrices, generalized Hadamard matrices and orthomorphism graphs of groups. In: \textit{JCMCC} 1 (1987), pp. 97–105.

\bibitem{GG2005} 
Gennian Ge. On $(g,4;1)$-difference matrices. In: \textit{Discrete Mathematics} 301.2–3 (2005), pp. 164 –174.

\bibitem{GGYM2010}
Gennian Ge, Ying Miao, and Xianwei Sun. Perfect difference families, perfect difference matrices, and related combinatorial structures. In: \textit{J. Combin. Designs} 18.6 (2010), pp. 415–449.

\bibitem{AGJS1979}
A. V. Geramita and J. Seberry. \textit{Orthogonal Designs: Quadratic forms and Hadamard matrices}. New York-Basel: Marcel Dekker, 1979. 460 pp.

\bibitem{JH1893}
J. Hadamard. Resolution d’une question relative aux determinants. In: \textit{Bull. des Sciences Mathematiques} 17 (1893), pp. 240–246

\bibitem{MHCL2010}
Masaaki Harada, Clement Lam, A. Munemasa, and Vladimir D. Tonchev. Classification of generalized Hadamard matrices $H(6,3)$ and quaternary Hermitian self-dual codes of length 18. In: \textit{ArXiv e-prints} (2010). Journal reference: Electronic J. Combin. 17 (2010), R171. arXiv: 1007.2555 [math.CO].

\bibitem{JH1997}
J.L. Hayden. Generalized Hadamard matrices. In: \textit{Des. Codes Cryptography} 12 (1997), pp. 69–73.

\bibitem{YH2014}
Yutaka Hiramine. On the non-existence of maximal difference matrices of deficiency 1. In: \textit{Designs, Codes and Cryptography} 72.3 (2014), pp. 627–635.

\bibitem{YHCS2013}
Yutaka Hiramine and Chihiro Suetake. On difference matrices of coset type. In: \textit{J. Combin. Theory} Ser. A 120.1 (2013), pp. 266–274.

\bibitem{KH2000}
K.J. Horadam. An introduction to cocyclic generalised Hadamard matrices. In: \textit{Discrete Applied Mathematics} 102.1-2 (2000), pp. 115 –131.

\bibitem{DJGG1986} 
Dieter Jungnickel and Gerhard Grams. Maximal difference matrices of order $\leq 10$. In: \textit{Discrete Mathematics} 58.2 (1986), pp. 199 –203.

\bibitem{DJ1980}
Dieter J. Jungnickel. On difference matrices and regular latin squares. In: \textit{Abh. Math. Sem. Hamburg} 50.1 (1980), pp. 219–231.

\bibitem{BK2009} 
B.R. Karlsson. Two-parameter complex Hadamard matrices for $N = 6$. In: \textit{J. Math. Phys}. 50 (2009), p. 82104.

\bibitem{BK2011}
B.R. Karlsson. Three-parameter complex Hadamard matrices of order 6. In: \textit{Linear Algebra Appl}. 434 (2011), p. 247.

\bibitem{PLFSPO2012}
Pekka H.J. Lampio, F. Sz\"{o}ll\H{o}si,and Ptric .R.J. \"{O}stergard. The quaternary complex Hadamard matrices of orders 10, 12, and 14. In: \textit{ArXiv e-prints} (2012). arXiv: 1204.5164 [math.CO].

\bibitem{PLPO2011}
Pekka H.J. Lampio and Patric R.J. \"{O}stergard. Classification of difference matrices over cyclic groups. In: \textit{J. Statistical Planning and Inference} 141.3 (2011), pp. 1194 –1207

\bibitem{WD1983}
Warwick de Launey. Generalised Hadamard matrices whose rows and columns form a group. In: \textit{Combinatorial Mathematics X}. Ed. by L.R.A. Casse. Vol. 1036. Lecture Notes in Mathematics. Berlin-Heidelberg-New York: Springer-Verlag, 1983, pp. 154–176.

\bibitem{1984b}
Warwick de Launey. On the non-existence of generalised Hadamard matrices. In: \textit{J. Statist. Plan. Inf}. 10 (1984), pp. 385–386.

\bibitem{WD1986a}
Warwick de Launey. A survey of generalised Hadamard matrices and difference matrices $D(k,\lambda;G)$ with large $k$. In: \textit{Utilitas Math}. 30 (1986), pp. 5–29.

\bibitem{WDL1986b}
Warwick de Launey. \textit{Some constructions for square GBRDs and some new infinite families of generalised Hadamard matrices}. Preprint. 1986.

\bibitem{WD1987b}
Warwick de Launey. On difference matrices, transversal designs, resolvable transversal designs and large sets of mutually orthogonal $F$-squares. In: \textit{J. Statist. Plan. Inf}. 16 (1987), pp. 107– 126.

\bibitem{WD1989}
Warwick de Launey. $GBRD$s: Some new constructions for difference matrices, generalised Hadamard matrices and balanced weighing matrices. In: \textit{Graphs and Combin}. 5 (1989), pp. 125–135.

\bibitem{WDL1992}
Warwick de Launey. Generalised Hadamard matrices which are developed modulo a group. In: \textit{Discrete Math}. 104 (1992), pp. 49–65.









\bibitem{WdeL07}
W. de Launey,  Generalized Hadamard matrices,   in C.J. Colbourn, J.H. Dinitz (Eds.), \textit{Handbook of Combinatorial Designs} (second ed.), Chapman \& Hall, CRC, Boca Raton, FL (2007), pp. 301–305

\bibitem{WDLED1992}
Warwick de Launey and J. E. Dawson. A note on the construction of $GH(4tq;EA(q))$ for $t = 1, 2$. In: \textit{Australas. J. Combin}. 6 (1992), pp. 177–186.

\bibitem{WDLED1994}
Warwick de Launey and J. E. Dawson. An asymptotic result on the existence of generalised Hadamard matrices. In: \textit{J. Combin. Theory}, A 65.1 (1994), pp. 158–163.

\bibitem{ZL1977}
Z.W. Liu. For the odd prime $p$ the construction of difference matrix for $OA(2p^2,2p+1;p,2)$. In: \textit{Acta Math. Appl. Sinica} 3 (1977), pp. 35–45.




\bibitem{La-Jolla}
La Jolla Difference Set Repository. 
URL \url{www.ccrwest.org/ds.html}. Viewed 2014:10:03.

\bibitem{MMFS2007}
M. Matolcsi and F. Sz\"{o}ll\H{o}si. Towards a classification of $6 \times 6$ complex Hadamard matrices. In: \textit{ArXiv Mathematics e-prints} (2007). arXiv: math/0702043.

\bibitem{JRS1979b}
Jennifer Seberry. Some remarks on generalized Hadamard matrices and theorems of Ra- jkundlia on SBIBDs. In: \textit{Combinatorial Mathematics VI}. Ed. by A.F. Horadam and W.D. Wallis. Vol. 748. Lecture Notes in Mathematics. Berlin-Heidelberg-New York: Springer-Verlag, 1979, pp. 154–164.

\bibitem{JRS1980a}
Jennifer Seberry. A construction for generalized Hadamard matrices. In: \textit{J. Statist. Plan. Inf.} 4 (1980), pp. 365–368.

\bibitem{JRS1986}
Jennifer Seberry. Generalized Hadamard matrices and colourable designs in the construction of regular GDDs with two and three association classes. In: \textit{J. Statistical Planning and Inference} 15.0 (1986-1987), pp. 237–246.

\bibitem{SebSTJOSEPH2015b}
Jennifer Seberry and N.A. Balonin,
Equivalence of the Existence of Hadamard Matrices and Cretan$(4t-1,2)$-Mersenne Matrices with Maximum Determinant, \textit{Journal of Theoretical and Computational Mathematics,} St. Joseph's College, Irinjalakuda, India,(submitted)


\bibitem{SFAWW}
Jennifer Seberry, Ken Finlayson, Sarah Spence Adams, Beata Wysocki, Tadeusz Wysocki and Tianbing Xia, The theory of quaternion orthogonal designs, In: \textit{IEEE Transactions on Signal Processing}, 56 1 (2008), 256-265. 

\bibitem{SY92}
 Jennifer Seberry and Mieko Yamada. Hadamard matrices, sequences, and block designs, \textit{Contemporary Design Theory: A Collection of Surveys}, J. H. Dinitz and D. R. Stinson, eds., John Wiley and Sons, Inc., 1992. pp. 431–560.

\bibitem{BSNS2010}
Bhu Dev Sharma and Norris Sookoo. Eigenvalues of the difference matrices of the Lee partition. In: \textit{J. Discrete Mathematical Sciences and Cryptography} 13.2 (2010), pp. 175–183.

\bibitem{SS1964}
S. Shrikhande. Generalized Hadamard matrices and orthogonal arrays strength two. In: \textit{Canad. J. Math.} 16 (1964), pp. 736–740.

\bibitem{DS1979}
D. Street. Generalised Hadamard matrices, orthogonal arrays and F-squares. In: \textit{Ars Combin}. 8 (1979), pp. 131–141.

\bibitem{FS2010b}
Ferenc Sz\"{o}ll\H{o}si. Exotic complex Hadamard matrices, and their equivalence. In: \textit{ArXiv e- prints} (2010). To appear in Cryptography and Communications: Discrete Structures, Boolean Functions and Sequences; arXiv: 1001.3062 [math.CO].

\bibitem{FS2011a}
Ferenc Sz\"{o}ll\H{o}si. Construction, classification and parametrization of complex Hadamard matrices. PhD thesis. The University of Wisconsin - Madison, 2011.

\bibitem{FS2012a} 
Ferenc Sz\"{o}ll\H{o}si. A note on the existence of $BH(19,6)$ matrices. In: \textit{ArXiv e-prints} (2012). arXiv: 1204.5166 [math.CO].

\bibitem{FS2012b}
Ferenc Sz\"{o}ll\H{o}si. On quaternary complex Hadamard matrices of small orders. In: \textit{ArXiv e-prints} (2012). Journal reference: Advances in Mathematics of Communications, 5 309–315 (2011). arXiv: 1204.5160 [math.CO].

\bibitem{DT2012}
Dobromir Todorov. Difference matrices toward MOLS of order 10. In: \textit{C. R. Acad. Bulgare Sci}. 65.8 (2012), pp. 1029–1034.

\bibitem{WSW}
W.D. Wallis, Anne Penfold Street and Jennifer Seberry Wallis, In: \textit{Combinatorics : Room Squares, Sum-free Sets, Hadamard Matrices}, 292, Lecture Notes in Mathematics, Springer--Verlag, Berlin--Heidelberg--New York, (1972), 508 pages. 

\bibitem{AW2000}
Arne Winterhof. On the non-existence of generalized Hadamard matrices. In: \textit{J. Statistical Planning and Inference} 84.1–2 (2000), pp. 337–342.

\bibitem{YZ2010} 
Yingshan Zhang. “Construction of difference matrices $D(r^m(r + 1), b^{m-1}a(b + 1); p)$. In: \textit{Advances and Applications in Discrete Mathematics} 6.2 (2010), pp. 77–100.

\bibitem{YZLD2002}
Yingshan Zhang, Lian Duan, Yiqiang Lu, and Zhongguo Zheng. Construction of generalized Hadamard matrices $D(r^m(r+1),r^m(r+1);p)$. In: \textit{J. Statistical Planning and Inference} 104.2 (2002), pp. 239 –258.

  \end{thebibliography}

\end{document}